# Strategic Adoption of 3D Printing in Multi-Product Supply Chains: Cost and Capacity Considerations

Mohammad E. Arbabian

Pamplin School of Business, University of Portland, Portland, Oregon, United States

arbabian@up.edu

# Strategic Adoption of 3D Printing in Multi-Product Supply Chains: Cost and Capacity Considerations

**Abstract**

This paper explores the integration of Additive Manufacturing (or 3D printing) into decentralized supply chains, focusing on the strategic decisions manufacturers and retailers make when facing capacity constraints. Using a Stackelberg game framework, we analyze how AM impacts traditional wholesale pricing across two scenarios: (1) a manufacturer producing two distinct products and (2) a case involving multiple products, each affected by AM's capacity limitations. For the first scenario, we derive sufficient conditions to find the equilibrium under a generic demand distribution, and for a uniform distribution, we fully derive the equilibrium and identify the critical cost threshold below which AM adoption is preferable. For the second scenario, we determine the conditions for the equilibrium, offering insights into the feasibility of AM relative to traditional manufacturing. Furthermore, we find that, contrary to some literature, even when AM incurs higher per-unit costs than traditional methods, it remains viable for some supply chains. However, limited AM capacity, particularly under high demand, may restrict its adoption, highlighting capacity constraints as a pivotal factor in AM decision-making. This study extends current research by examining multi-product and capacity-driven scenarios, providing valuable guidance for supply chain managers weighing the benefits of traditional manufacturing against those of AM.

## 1. Introduction

The advent of 3D printing, or additive manufacturing, has revolutionized manufacturing, which traditionally utilizes production processes such as molding, casting, and machining. This has enabled firms to transition from large-scale, standardized manufacturing to highly customizable and flexible production. This shift is especially important in supply chain management, where manufacturers and retailers face increasing pressures to meet dynamic consumer demand, reduce lead times, and minimize costs. Traditional manufacturing processes, which rely heavily on economies of scale, often require large production runs to be cost-effective. However, 3D printing offers a solution by allowing on-demand production of smaller batches with minimal setup costs. As the technology advances, many industries, from aerospace[1] to healthcare[2], are adopting 3D printing to streamline their supply chains and respond more efficiently to distinct market demands.

---
[1] https://formlabs.com/blog/additive-manufacturing-3d-printing-in-aerospace/

In a traditional supply chain setup, the manufacturer sells products to the retailer at a wholesale price, and the retailer, in turn, determines the optimal order quantities based on customer demand. In other words, the manufacturer sets the wholesale price, and the retailer manages the uncertainties of customer demand through inventory and pricing strategies. This type of supply chain typically operates under a wholesale price contract, which is a common arrangement in decentralized supply chains. In this scenario, the manufacturer and retailer are independent entities with separate objectives: the manufacturer aims to maximize profit by setting the wholesale price, while the retailer seeks to optimize order quantities to balance the cost of purchasing inventory with the uncertain demand from customers (Lariviere and Porteus, 2001). In this study, we build on this traditional framework by introducing 3D printing as an option for the manufacturer.

This study is motivated by the challenges faced by long-established companies that produce spare parts for legacy products. To meet demand for these parts, such companies often maintain outdated production lines, which can be costly and inefficient. By adopting 3D printing, however, these companies could close these older production lines while still manufacturing spare parts for multiple legacy products. A notable example is BMW, which has utilized 3D printing to replace traditional production lines for older vehicle models, reducing operational costs and enhancing production flexibility[3]. Nike also reported 13% growth in revenue in the last quarter of 2014 because of adopting 3D printing. "Nike made 3D printed cleats for the 2014 Super Bowl". Hasboro, also, recently adopted 3D printing (partnered with 3D Systems), and its CEO stated that "We believe 3D printing offers endless potential to bring incredible new play experiences for kid" [4].

Therefore, to utilize 3D printing, the manufacturer must first decide whether to adopt 3D printing or not. If 3D printing is not adopted, the problem reduces to the traditional wholesale price contract, where the manufacturer sets the wholesale price and the retailer orders products accordingly. However, if 3D printing is adopted, the manufacturer must decide at what wholesale price to sell the products to the retailer, considering both the benefits of customization and the limitations imposed by production capacity. For instance, if BMW chose not to implement 3D printing, it would need to maintain the legacy car production line, and dealerships would rely on traditional wholesale pricing contracts to determine optimal order quantities from BMW. Alternatively, if BMW adopts 3D printing technology, it must establish appropriate pricing for spare parts produced using this method for legacy cars.

---

[2] https://www.synopsys.com/glossary/what-is-medical-3d-printing.html
[3] https://www.automationworld.com/factory/3d-printing-additive-manufacturing/article/33018717/bmw-uses-3d-printing-to-manufacture-powertrain-components
[4] https://www.techrepublic.com/article/3d-printing-10-companies-using-it-in-ground-breaking-ways/

The integration of 3D printing into a manufacturer-retailer supply chain introduces several complexities. On the one hand, 3D printing allows for the production of multiple products, which can be particularly advantageous in markets where demand is volatile or difficult to predict. On the other hand, the high fixed costs of 3D printing equipment and the potential for capacity constraints require the manufacturer to carefully evaluate whether the benefits of 3D printing outweigh its costs. The tradeoff we examine is between the flexibility of 3D printing and its capacity constraints. That is, while 3D printing offers the flexibility to produce various unique products without incurring setup costs, it operates much more slowly than traditional manufacturing methods, thus limiting its overall capacity.

This paper considers two distinct scenarios to analyze these complexities and tradeoffs. In the first scenario, we focus on a supply chain with *two distinct products*, which allows for a detailed analysis of the Stackelberg equilibrium. The Stackelberg framework is particularly relevant for this analysis because it captures the hierarchical decision-making structure between the manufacturer (leader) and the retailer (follower). The manufacturer, acting as the leader, 1) decides whether 3D printing technology should be adopted or not, and 2) sets the wholesale price for each product. The retailer, acting as the follower, observes these prices and determines the optimal order quantities based on expected customer demand. In the first scenario, we explore how the manufacturer's decision to adopt 3D printing affects the equilibrium outcomes, and we identify the sufficient conditions under which 3D printing is preferred over traditional manufacturing. Furthermore, for a uniform distribution, we fully derive the equilibrium. In the second scenario, we extend the model to consider the production of *multiple products*, where we are able to derive the sufficient condition to find the equilibrium. In both scenarios (i.e., two product and multi-product setting), the manufacturer must decide not only whether to adopt 3D printing but also how to allocate production capacity across different products. Capacity constraints become a critical factor in these scenarios, as 3D printers typically have limited production capacity, and the manufacturer must balance the production of multiple products to maximize profit. As the result, the retailer's ordering decisions are influenced not only by the wholesale price but also by the availability of different products, which may be subject to capacity limitations.

Our model incorporates several key factors that influence the adoption of 3D printing, including the fixed costs associated with purchasing 3D printers, the per-unit production costs of both 3D printing and traditional manufacturing, and the stochastic nature of customer demand. The decision-making process in this decentralized supply chain is further complicated by the fact that the manufacturer and retailer have independent objectives, with the manufacturer seeking to maximize profit by setting the optimal wholesale price and the retailer seeking to optimize order quantities to meet uncertain demand. By comparing the equilibrium outcomes in both scenarios, we provide insights into the conditions under which 3D printing is a viable option for the manufacturer, and how the retailer's ordering strategy is influenced by the manufacturer's decision.

Arbabian (2022) examined a single-product version of this problem and argued that 3D printing is not viable when the variable costs of 3D printing exceed those of traditional manufacturing. However, this study demonstrates that 3D printing can still be beneficial in certain cases, even when the variable costs are higher. Specifically, we show that for multi-product settings, the advantages of flexibility and the ability to offset costs across multiple products can make 3D printing profitable, *despite higher per-unit costs*.

This study contributes to the existing literature on 3D printing in supply chains in several important ways. First, we extend the analysis of 3D printing to a decentralized supply chain model, providing insights into how the manufacturer's decision to adopt 3D printing affects the retailer's behavior. While much of the existing research on 3D printing focuses on its benefits in isolated firms or single-stage production processes, this paper examines the broader implications of 3D printing in a multi-stage supply chain. Second, we introduce a multi-product model that incorporates capacity constraints and product-specific costs, offering a more comprehensive analysis of how 3D printing affects supply chain dynamics when multiple products are involved. Finally, we provide practical insights for supply chain managers, highlighting the trade-offs between traditional manufacturing and 3D printing in terms of cost efficiency, capacity utilization, and responsiveness to customer demand.

The remainder of this paper is structured as follows. In Section 2, we review the relevant literature on 3D printing in supply chains and decentralized decision-making. Section 3 introduces the benchmark for our proposed models. Section 4 lays out our theoretical model for the two-product scenario, providing a detailed analysis of the Stackelberg equilibrium. In Section 5, we extend the model to account for multiple products and examine how capacity constraints and product-specific costs influence the adoption of 3D printing in the supply chain. Section 6 presents numerical studies that illustrate the impact of key parameters on the equilibrium outcomes with managerial insights and practical recommendations for supply chain managers considering the integration of 3D printing into their operations.

## 2. Literature Review

The intersection of 3D printing (additive manufacturing) and supply chain management has garnered increasing attention in recent years due to the transformative potential of this technology. This section reviews the relevant literature in three domains: 3D printing in operations management, supply chain contracting, and supply chain investment, positioning this study within these streams.

### 2.1 3D Printing in Operations Management

The literature on 3D printing's impact on supply chain dynamics is still emerging, with a particular focus on its ability to reduce lead times, increase customization, and lower inventory costs. Conner et al. (2014) presented a framework that examines where 3D printing is economically feasible, identifying that its advantages are most pronounced in cases with high product complexity, high customization, or low production volume. Dong et al. (2017) compared 3D printing with traditional flexible manufacturing technologies, demonstrating how 3D printing can offer full flexibility at lower costs. Their findings suggest that while traditional flexible technologies are often more expensive to scale, 3D printing offers a cost-effective alternative. Westerweel et al. (2018) explored how 3D printing can reduce inventory and lead times for spare parts in supply chains. Their research emphasizes the value of 3D printing in situations where maintaining large inventories of parts is costly or infeasible. Similarly, Song and Zhang (2019) examined the shift from make-to-stock to make-to-order production strategies enabled by 3D printing. They found that firms adopting 3D printing can move toward more responsive production models, which reduce excess inventory. Chen et al. (2020) expanded on this by studying how 3D printing can be integrated across multiple channels in a supply chain, particularly in e-commerce. They highlighted the flexibility of 3D printing to meet diverse demand from both online and offline channels. Sethuraman et al. (2018) focused on the end-user adoption of 3D printing technology, or "personal fabrication," where consumers print products themselves. Their research provides insights into how 3D printing disrupts traditional supply chain models by empowering consumers to bypass retailers entirely.

The main difference between this stream of literature and our study is that we mainly focus on a decentralized supply chain while all these studies focus on a centralized supply chain.

Perhaps, one of the most relevant studies to this paper is Arbabian and Wagner (2020), who examined a single-product supply chain with 3D printing. They highlighted key economic conditions that favor the adoption of 3D printing by manufacturers. Their findings show that, in some cases, retailer adoption of 3D printing can eliminate double marginalization, a result echoed in our study. Several other studies also contribute to this growing body of research. For example, Hall (2016) examined the economic conditions under which manufacturers can use 3D printers to produce other 3D printers, a concept known as self-replication. This idea aligns with broader research on 3D printing's ability to enable decentralized manufacturing (Kretschmer, 2015), where production can occur closer to the point of demand, reducing transportation costs and increasing responsiveness.

while these papers focus on producing and selling only one type of product, we focus is on supply chains where more than one type of product is produced and sold.

Furthermore, Ahmed, Heese, and Kay (2023) demonstrate how stochastic optimization can optimize AM deployment in manufacturing networks, addressing uncertainties in demand and production costs. Similarly, Cantini et al. (2022) propose a decision support system that aids organizations in deciding when to use AM for spare parts, balancing cost, lead time, and service levels. Both studies underline AM's potential to improve supply chain efficiency and responsiveness to varying demands. The COVID-19 pandemic has further influenced supply chain strategies, as Ivanov (2021) introduces the concept of "supply chain viability," advocating for adaptability through redundancy, flexibility, digitalization, and sustainability. Complementing this, Sgarbossa et al. (2021) analyze the comparative effectiveness of conventional versus AM for spare parts management, finding that AM offers notable advantages for low-volume, high-variability items. Additionally, Top et al. (2023) explore AM's environmental benefits and operational challenges in promoting sustainable manufacturing, concluding that while there are obstacles to adoption, AM significantly contributes to reducing waste and enhancing production flexibility.

The main difference between this body of literature and our research is the cost structure. That is, in our paper, we utilize a wholesale price contract where the cost of 3D printed products is endogenous.

Finally, Kucukkoc (2019) presents Mixed-Integer Linear Programming (MILP) models designed to minimize makespan in additive manufacturing (AM) machine scheduling, focusing on optimizing job sequences to enhance scheduling efficiency and balance workloads across machines. This optimization can significantly improve operational efficiency in AM environments. Meanwhile, Peron et al. (2022) highlight the role of AM in bolstering supply chain resilience, particularly during the COVID-19 pandemic, by enabling localized production and reducing reliance on global suppliers. They argue that AM allows companies to quickly adapt to disruptions by facilitating on-demand production of critical components, ultimately promoting a sustainable and flexible approach to managing future supply chain challenges.

While these papers study 3D printing in the context of supply chain, their basic assumption is that 3D printing is already adopted, and they aim to improve the efficiency of a supply chain where 3D printing is utilized. In our paper, however, adopting 3D printing technology is a decision variable.

## 2.2. Supply Chain Contracting

Supply chain contracting literature has long focused on how contracts can coordinate decentralized supply chains, often with the goal of eliminating the double marginalization effect. Pasternack (1985) first showed that buy-back contracts can coordinate supply chains by aligning the incentives of manufacturers and retailers. Cachon and Lariviere (2005) extended this by demonstrating that revenue-sharing contracts could achieve similar coordination outcomes. Tsay (1999) added that quantity-flexibility contracts can also coordinate supply chains effectively. However, 3D printing introduces new complexities into this traditional contracting literature. Arbabian (2022) argued that if the variable cost of 3D printing exceeds that of traditional manufacturing, 3D printing should not be adopted. Our study, however, challenges this notion by showing that even when the cost of 3D printing is higher, the benefits of customization and flexibility can justify its adoption in multi-product supply chains.

More recent studies have examined how 3D printing can be integrated into existing supply chain contracts. For instance, Chen et al. (2019) explored how contracts can be structured to incentivize the adoption of 3D printing in decentralized supply chains. Their research highlights the importance of wholesale-price contracts in coordinating supply chains that include both traditional manufacturing and 3D printing.

## 2.3 Supply Chain Investment

The adoption of 3D printing is analogous to investing in manufacturing technologies that reduce production costs or increase capacity. This literature dates back to Porteus (1985), who studied the optimization of investment in setup costs within an EOQ model. Fine and Freund (1990) extended this analysis to flexible manufacturing systems, while Van Mieghem (1998) explored optimal investment strategies in flexible manufacturing environments. Goyal and Netessine (2007) investigated the interaction between competition and investment in technology, finding that firms often under-invest in technology due to competitive pressures.

Our study fits within this stream by analyzing the conditions under which firms should invest in 3D printing technology. Specifically, we model the fixed cost of acquiring a 3D printer and examine how firms should make investment decisions based on demand variability, product complexity, and cost structures. Similar to Ge et al. (2014), who examined investment strategies for suppliers and manufacturers, we consider how investment in 3D printing affects both upstream and downstream supply chain members.
Finally, Bernstein and Kök (2009) analyzed how different contracts influence investment decisions in supply chains with multiple suppliers. Their findings suggest that the allocation of risk between suppliers and manufacturers can significantly impact investment in new technologies. Our research extends this by showing how 3D printing investments can be coordinated between manufacturers and retailers through wholesale-price contracts.

### 3. Benchmark Model Without 3D Printing

In this section, we introduce our benchmark model, reflecting the traditional scenario in which firms lack access to 3D printing (i.e., Scenario b). This model serves two main objectives: 1) We will showcase our findings for scenarios where firms do possess 3D printing capabilities, as outlined in this section. 2) We will conduct a comparative analysis of firm outcomes between the benchmark model and those incorporating 3D printing.

In our benchmark model, a scenario involving a solitary manufacturer selling a single product to a sole retailer via a wholesale-price contract is studied. The unit wholesale price is denoted by $w \geq 0$. This setup is rooted in the model presented by Lariviere and Porteus (2001). Conversely, the retailer distributes the products to a continuous stochastic demand $D$ characterized by density function $f$ and cumulative distribution function $F$, with a range on $[0, U]$. The manufacturing variable cost is denoted by $c_m > 0$, and the retailer's unit selling price is shown by $r > 0$. To prevent trivial solutions, we assume that $c_m \leq w \leq r$. The retailer's single decision, contingent upon $w$, is to ascertain the order quantity $q(w) \geq 0$ that maximizes his expected profit.

$$\pi_R^b(w) = \max_{q \geq 0} E[r \min\{q, D\}] - wq. \tag{1}$$

This model is solved by the famous classic Newsvendor solution $q_{nv}(w) = F^{-1}\left(1 - \frac{w}{r}\right)$. $\pi_R^b(w)$ represents the *retailer's* profit for our *benchmark*.

In this Stackelberg game, manufacturer's decision is setting the appropriate wholesale price $w$ that maximizes her profit. Manufacturer's Profit is

$$\pi_M^b = \max_{w \geq 0} (w - c_m) q_{nv}(w). \tag{2}$$

To solve this problem, we make the following assumption similar to (Lariviere and Porteus 2001).

**Assumption 1**. The demand distribution $F$ satisfies the increasing generalized failure rate constraint.

Generalized failure rate for the demand distribution $F$ is defined as $\frac{xf(x)}{1-F(x)}$, which Assumption 1 indicates that is increasing with respect to $x$. Most common distributions like uniform, gamma, Pareto fulfill this constraint. For a more comprehensive list, see Banciu and Mirchandani (2013).

**Proposition 1**. *(Lariviere and Porteus 2001) The manufacturer's profit function $(w - c_m)q_{nv}(w)$ is strictly unimodal in w. The optimality condition is:*

$$(1 - F(q))\left(1 - \frac{qf(q)}{1 - F(q)}\right) = \frac{c_m}{r} \tag{3}$$

We define $q^0 = q_{nv}(c_m)$ as the unique solution to Equation (3). Therefore, $w^0 = r(1 - F(q^0))$ denotes the associated optimal wholesale price.

## 4. Manufacturer May Purchase a 3D Printer—Two Products

In this section, we investigate a scenario where the manufacturer uses traditional manufacturing methods to produce two distinct products. The unit production cost for the two products is $(c_{m,1}, c_{m,2})$. The manufacturer sells the products to the retailer at the wholesale price of $(w_1, w_2)$. The retailer, in turn, faces stochastic demand for each product with density functions of $(f_1(.), f_2(.))$ and cumulative distribution functions of $(F_1(.), F_2(.))$. The retailer sells the products to the final customer at the price of $(r_1, r_2)$. In this case, following our benchmark case, the retailer's problem is:

$$\pi_R^b(w) = \max_{q_i \geq 0} E[r_1 \min\{q_1, D_1\}] - w_1 q_1 + E[r_2 \min\{q_2, D_2\}] - w_2 q_2.$$

Following Proposition 1, the optimal solution to the retailer's problem is $(q_{nv}(w_1), q_{nv}(w_1)) = \left(F_1^{-1}\left(1 - \frac{w_1}{r_1}\right), F_2^{-1}\left(1 - \frac{w_2}{r_2}\right)\right)$.

Next, in presence of 3D printing, the manufacturer has the option of purchasing a 3D printer with a fixed cost $K > 0$ and unit printing costs $(c_{p,1}, c_{p,2})$, to replace her traditional manufacturing. This is typically observed in situations characterized by high complexity, significant customizability, and/or low production volumes. An excellent illustration of this concept is evident in the actions of car manufacturers like BMW, which are integrating 3D printing technology into their production processes. Usually, even after a particular car model ceases production, many manufacturers retain their old production lines to manufacture spare parts for these discontinued models. Transitioning to 3D printing allows companies to phase out various outdated car production lines and replace them with 3D printing technology.

We let $v \in \{0,1\}$ denote the manufacturer's binary decision whether to adopt 3D printing technology. If 3D printing is adopted ($v = 1$), there will be no traditional manufacturing (i.e., it is illogical for BMW to retain its old car production lines solely for the purpose of manufacturing spare parts while concurrently implementing 3D printing technology for the same purpose). One downside of 3D printing is that its capacity $Q$ is limited compared to traditional manufacturing techniques. Alternatively, if 3D printing is not adopted ($v = 0$), then only traditional manufacturing would be possible. Therefore, in presence of 3D printing, the manufacturer's profit-maximizing problem is

$$\pi_M^{3D} = \max_{w_i \geq 0, v \in \{0,1\}} \left(w_1 - vc_{p,1} - (1-v)c_{m,1}\right)q_{nv}(w_1) + \left(w_2 - vc_{p,2} - (1-v)c_{m,2}\right)q_{nv}(w_2) - Kv \quad (4)$$

$$s.t. \quad q_{nv}(w_1) + q_{nv}(w_2) \leq (1-v)M + vQ.$$

The first term in the objective function in (4) represents the profit from selling product 1. If 3D printing is *not* adopted ($v = 0$), this term simplifies to that of the benchmark for product 1. If 3D printing, however, is opted for ($v = 1$), then the first term simplifies to $(w_1 - c_{m,1})q_{nv}(w_1)$, which is manufacturer's profit from selling 3D printed product 1s. Similar analogy holds for the second term. The last term is the fixed cost of adopting/investing in 3D printing technology. In the constraint in (4), $M$ represents a large numerical value. This constraint guarantees that should 3D printing be adopted, the cumulative optimal order quantities determined by the retailer adheres to the capacity limit imposed by 3D printing. In the subsequent sections, we delve into problem (4) to derive the Stackelberg equilibrium.

### 4.1. 3D Printing Is Not Adopted

In this section, we study the scenario where 3D printing is not adopted (i.e., $v = 0$), in which case, the problem in (4) simplifies to

$$\pi_M^{3D} = \max_{w_i \geq 0} \quad (w_1 - c_{m,1})q_{nv}(w_1) + (w_2 - c_{m,2})q_{nv}(w_2) \tag{5}$$
$$s.t. \quad q_{nv}(w_1) + q_{nv}(w_2) \leq M,$$

Note that, in this case, the constraint is redundant. This problem is similar to that of the benchmark. The only difference is that the retailer orders two district products. Therefore, the optimal solution can be achieved from the following Proposition.

**Proposition 2**. *The unique optimal solution to the manufacturer's Problem (5) is*
$$(q_1^0, q_2^0) \triangleq \left(q_b(c_{m,1}), q_b(c_{m,2})\right) = \left(F_1^{-1}\left(1 - \frac{c_{m,1}}{r_1}\right), F_2^{-1}\left(1 - \frac{c_{m,2}}{r_2}\right)\right)$$
$$(w_1^0, w_2^0) \triangleq \left(r_1\left(1 - F_1(q_1^0)\right), r_2\left(1 - F_2(q_2^0)\right)\right)$$

As one may observe, because in this case 3D printing is not adopted, the optimal order quantities and the optimal wholesale prices are similar to that of the benchmark (i.e., Proposition 1). To gain deeper insight into the implications of Proposition 2, we distill its findings by applying them to a uniform distribution. Moreover, leveraging a uniform distribution aids in extracting further insights from the results presented in the subsequent sections.

**Example 1.** *To explore the implications of the problem in (4), in this section, we advance our analysis focusing on the scenario where $F_1 \sim Uniform\ [0, U_1]$, $F_2 \sim Uniform\ [0, U_2]$.*
$$(q_1^0, q_2^0) = \left(-\frac{U_1(c_{m,1} - r_1)}{2r_1}, -\frac{(c_{m,2} - r_2)U_2}{2r_2}\right),$$
$$(w_1^0, w_2^0) = \left(r_1\left(1 + \frac{c_{m,1} - r_1}{2r_1}\right), r_2\left(1 + \frac{c_{m,2} - r_2}{2r_2}\right)\right),$$
$$\pi_M^0 = \frac{U_1 r_1^2\ r_2 + (U_2 r_2^2 + (-2U_1 c_{m,1} - 2U_2 c_{m,2})r_2 + U_2 c_{m,2}^2)r_1 + U_1 c_{m,1}^2 r_2}{4r_1 r_2}.$$

### 4.2. 3D Printing Is Adopted

In this section, we study the scenario where 3D printing is adopted (i.e., $v = 1$). In this case, the problem in (4) simplifies to

$$\pi_M^{3D} = \max_{w_i \geq 0} \quad (w_1 - c_{p,1})q_{nv}(w_1) + (w_2 - c_{p,2})q_{nv}(w_2) - K \tag{6}$$
$$s.t. \quad q_{nv}(w_1) + q_{nv}(w_2) \leq Q,$$

Given the constraint in (6), we divide the solution space into two regions.

**Region 1.** $Q > q_{nv}(w_1) + q_{nv}(w_2)$. In this case, the constraint in (6) is redundant, and there is an interior solution to the problem in (6). This case represents a scenario where 3D printing is utilized but *not* at its maximum capacity. The following lemma maximizes manufacture's profit in this region.

**Lemma 1.** *The unique optimal solution to the manufacturer's Problem (6) is*

$$(q_1^1, q_2^1) \triangleq (q_b(c_{p,1}), q_b(c_{p,2})) = \left(F_1^{-1}\left(1 - \frac{c_{p,1}}{r_1}\right), F_2^{-1}\left(1 - \frac{c_{p,2}}{r_2}\right)\right)$$

$$(w_1^1, w_2^1) \triangleq \left(r_1\left(1 - F_1(q_1^1)\right), r_2\left(1 - F_2(q_2^1)\right)\right)$$

$$\pi_M^1 \triangleq (w_1^1 - c_{p,1})q_b(c_{p,1}) + (w_2^1 - c_{p,2})q_b(c_{p,2}) - K$$

To gain deeper insight into the implications of Lemma 1, we distill its findings by applying them to a uniform distribution.

**Example 2.** *Following Example 1, we simplify our results in Lemma 1 for $F_1 \sim Uniform\ [0, U_1]$, $F_2 \sim Uniform\ [0, U_2]$.*

$$(q_1^1, q_2^1) \triangleq \left(\frac{U_1}{2}\left(1 - \frac{c_{p,1}}{r_1}\right), \frac{U_2}{2}\left(1 - \frac{c_{p,2}}{r_2}\right)\right),$$

$$(w_1^1, w_2^1) \triangleq \left(r_1\left(1 + \frac{c_{p,1} - r_1}{2r_1}\right), r_2\left(1 + \frac{c_{p,2} - r_2}{2r_2}\right)\right),$$

$$\pi_M^1 \triangleq \frac{U_1 r_1^2 r_2 + (r_2^2 U_2 + (-2U_1 c_{p,1} - 2U_2 c_{p,2})r_2 + U_2 c_{p,2}^2)r_1 + U_1 c_{p,1}^2 r_2}{4 r_1 r_2} - K.$$

As anticipated, the outcomes here mirror those of Example 1, albeit with the substitution of $c_{m,i}$ by $c_{p,i}$.

**Region 2.** $q_{nv}(w_1) + q_{nv}(w_2) \geq Q$. In this case, the constraint in (6) is binding, and there is a corner solution to the problem in (6). This case represents a scenario where 3D printing *is* utilized at its maximum capacity. Therefore, $q_{nv}(w_1) + q_{nv}(w_2) = Q \rightarrow q_{nv}(w_2) = Q - q_{nv}(w_1)$. Therefore, the problem in (6), simplifies to

$$\pi_M^{3D} = \max_{w_i \geq 0}\ (w_1 - c_{p,1})q_{nv}(w_1) + (w_2 - c_{p,2})(Q - q_{nv}(w_1)) - K \tag{7}$$

Next, the retailer's optimal order quantity is $q_{nv}(w_2) = F_2^{-1}\left(1 - \frac{w_2}{r_2}\right) = Q - q_{nv}(w_1)$. Therefore, $w_2 = r_2\left(1 - F_2(Q - q_{nv}(w_1))\right)$. Finally, the problem in (7), simplifies to

$$\pi_M^{3D} = \max_{q_1} \ (r_1(1-F_1(q_1)) - c_{p,1})q_1 + (r_2(1-F_2(Q-q_1)) - c_{p,2})(Q-q_1) - K \quad (8)$$

**Lemma 2**. *The optimal solution to the manufacturer's Problem (8) is the unique solution, $q_1^{1Q}$, to the following equation.*

$$r_1(1-F_1(q_1))\left[1 - \frac{f_1(q_1)q_1}{1-F_1(q_1)}\right] - r_2(1-F_2(Q-q_1))\left[1 - \frac{f_2(Q-q_1)(Q-q_1)}{1-F_2(Q-q_1)}\right] = c_{p,1} - c_{p,2} \quad (9)$$

*And, $q_2^{1Q} = Q - q_1^{1Q}$. Furthermore, $w_i^{1Q} = r_i\left(1 - F_i(q_i^{1Q})\right)$.*

To gain deeper insight into the implications of Lemma 2, we distill its findings by applying them to a uniform distribution.

**Example 3.** *Following Example 1, we simplify our results in Lemma 2 for $F_1 \sim Uniform\ [0, U_1]$, $F_2 \sim Uniform\ [0, U_2]$.*

$$(q_1^{1Q}, q_2^{1Q}) = \left(\frac{U_1(2Qr_2 - U_2 c_{p,1} + U_2 c_{p,2} + U_2 r_1 - U_2 r_2)}{2(U_1 r_2 + U_2 r_1)}, \frac{\left((c_{p,1} - c_{p,2} - r_1 + r_2)U_1 + 2r_1 Q\right)U_2}{2(U_1 r_2 + U_2 r_1)}\right),$$

$$(w_1^{1Q}, w_2^{1Q}) = \left(-\frac{r_1\left((-c_{p,1} + c_{p,2} - r_1 - r_2)U_2 + 2r_2(Q - U_1)\right)}{2(U_1 r_2 + U_2 r_1)}, -\frac{r_2\left((c_{p,1} - c_{p,2} - r_1 - r_2)U_1 + 2r_1(Q - U_2)\right)}{2(U_1 r_2 + U_2 r_1)}\right),$$

$$\pi_M^{1Q} = \frac{\left((c_{p,1} - c_{p,2} - r_1 + r_2)^2 U_2 - 4Qr_2(c_{p,1} - r_1)\right)U_1 - 4Qr_1\left((c_{p,2} - r_2)U_2 + Qr_2\right)}{4U_1 r_2 + 4U_2 r_1} - K.$$

Next, the optimal solution to Problem (6) is the maximum of the solutions in Region 1 and Region 2, which is derived in Proposition 3.

**Proposition 3**. *The unique optimal solution to the manufacturer's Problem (6) is*

$$(q_1, q_2) = \begin{cases} (q_1^{1Q}, q_2^{1Q}), & q_b(c_{p,1}) + q_b(c_{p,2}) > Q \\ (q_b(c_{p,1}), q_b(c_{p,2})), & \text{otherwise.} \end{cases},$$

*and*

$$(w_1^{1Q}, w_2^{1Q}) = \begin{cases} \left(r_1\left(1 - F_1(q_1^{1Q})\right), r_2\left(1 - F_2(q_2^{1Q})\right)\right), & q_b(c_{p,1}) + q_b(c_{p,2}) > Q \\ \left(r_1\left(1 - F_1(q_b(c_{p,1}))\right), r_2\left(1 - F_2(q_b(c_{p,2}))\right)\right), & \text{otherwise.} \end{cases}$$

**Example 4.** *In this example, we simplify the results in Proposition 3 for $F_1 \sim Uniform\ [0, U_1]$, $F_2 \sim Uniform\ [0, U_2]$.*

$(q_1, q_2)$

$$= \begin{cases} \left(\dfrac{U_1}{2}\left(1-\dfrac{c_{p,1}}{r_1}\right), \dfrac{U_2}{2}\left(1-\dfrac{c_{p,2}}{r_2}\right)\right), & \dfrac{U_1}{2}\left(1-\dfrac{c_{p,1}}{r_1}\right)+\dfrac{U_2}{2}\left(1-\dfrac{c_{p,2}}{r_2}\right) < Q \\ \left(\dfrac{U_1(2Qr_2 - U_2 c_{p,1} + U_2 c_{p,2} + U_2 r_1 - U_2 r_2)}{2(U_1 r_2 + U_2 r_1)}, \dfrac{\left((c_{p,1} - c_{p,2} - r_1 + r_2)U_1 + 2r_1 Q\right)U_2}{2(U_1 r_2 + U_2 r_1)}\right), & \text{otherwise.} \end{cases}$$

$(w_1, w_2)$

$$= \begin{cases} \left(r_1\left(1+\dfrac{c_{p,1}-r_1}{2r_1}\right), r_2\left(1+\dfrac{c_{p,2}-r_2}{2r_2}\right)\right), & \dfrac{U_1}{2}\left(1-\dfrac{c_{p,1}}{r_1}\right)+\dfrac{U_2}{2}\left(1-\dfrac{c_{p,2}}{r_2}\right) < Q \\ \left(-\dfrac{r_1\left((-c_{p,1}+c_{p,2}-r_1-r_2)U_2+2r_2(Q-U_1)\right)}{2(U_1 r_2 + U_2 r_1)}, -\dfrac{r_2\left((c_{p,1}-c_{p,2}-r_1-r_2)U_1+2r_1(Q-U_2)\right)}{2(U_1 r_2 + U_2 r_1)}\right), & \text{otherwise.} \end{cases}$$

$\pi_M^{3D}$

$$= \begin{cases} \dfrac{U_1 r_1^2 r_2 + \left(r_2^2 U_2 + (-2U_1 c_{p,1} - 2U_2 c_{p,2})r_2 + U_2 c_{p,2}^2\right)r_1 + U_1 c_{p,1}^2 r_2}{4 r_1 r_2} - K, & \dfrac{U_1}{2}\left(1-\dfrac{c_{p,1}}{r_1}\right)+\dfrac{U_2}{2}\left(1-\dfrac{c_{p,2}}{r_2}\right) < Q \\ \dfrac{\left((c_{p,1}-c_{p,2}-r_1+r_2)^2 U_2 - 4Qr_2(c_{p,1}-r_1)\right)U_1 - 4Qr_1\left((c_{p,2}-r_2)U_2 + Qr_2\right)}{4U_1 r_2 + 4U_2 r_1} - K, & \text{otherwise.} \end{cases}$$

Following Example 4, one can observe that the optimal order quantities are decreasing in the unit 3D printing cost (i.e. $\dfrac{\partial q_i}{\partial c_{p,i}} < 0$), and the optimal wholesale prices are increasing in the unit 3D printing cost (i.e. $\dfrac{\partial w_i}{\partial c_{p,i}} > 0$). This implies that the manufacturer should increase the wholesale price as 3D printing becomes more expensive. At the same time, the retailer should decrease the order quantity because the product becomes more expensive to purchase.

### 4.3. The Equilibrium

The optimal solution to the manufacturer's Problem (4) is the maximum of Problem (5) and Problem (6), which is derived in the following proposition.

**Proposition 4**. *The equilibrium to Problem (4) is:*

$$(q_1^*, q_2^*, v^*) = \begin{cases} (q_1^{1Q}, q_2^{1Q}, 1), & q_b(c_{p,1}) + q_b(c_{p,2}) \geq Q \text{ and } \pi_M^{1Q} > \pi_M^0 \\ (q_b(c_{p,1}), q_b(c_{p,2}), 1), & q_b(c_{p,1}) + q_b(c_{p,2}) < Q \text{ and } \pi_M^1 > \pi_M^0 \\ (q_b(c_{m,1}), q_b(c_{m,2}), 0), & \text{otherwise} \end{cases}$$

$$w_i^* = r_i(1 - F_i(q_i^*))$$

To have a better understanding of the results on Proposition 4, we simplify them for a uniform distribution in the following example.

**Example 5**: *For a Uniform distribution, the Equilibrium and the optimal profit are as follows:*

$$(q_1^*, q_2^*, v^*)$$

$$= \begin{cases} \left(\dfrac{U_1(2Qr_2 - U_2c_{p,1} + U_2c_{p,2} + U_2r_1 - U_2r_2)}{2(U_1r_2 + U_2r_1)}, \dfrac{\left((c_{p,1} - c_{p,2} - r_1 + r_2)U_1 + 2r_1Q\right)U_2}{2(U_1r_2 + U_2r_1)}, 1\right), & \dfrac{U_1}{2}\left(1 - \dfrac{c_{p,1}}{r_1}\right) + \dfrac{U_2}{2}\left(1 - \dfrac{c_{p,2}}{r_2}\right) > Q \\ \text{and } \dfrac{\left((c_{p,1} - c_{p,2} - r_1 + r_2)^2 U_2 - 4Qr_2(c_{p,1} - r_1)\right)U_1 - 4Qr_1\left((c_{p,2} - r_2)U_2 + Qr_2\right)}{4U_1r_2 + 4U_2r_1} - K > \dfrac{U_1r_1^2 r_2 + (U_2r_2^2 + (-2U_1c_{m,1} - 2U_2c_{m,2})r_2 + U_2c_{m,2}^2)r_1 + U_1c_{m,1}^2 r_2}{4r_1 r_2} \\ \left(\dfrac{U_1}{2}\left(1 - \dfrac{c_{p,1}}{r_1}\right), \dfrac{U_2}{2}\left(1 - \dfrac{c_{p,2}}{r_2}\right), 1\right), & \dfrac{U_1}{2}\left(1 - \dfrac{c_{p,1}}{r_1}\right) + \dfrac{U_2}{2}\left(1 - \dfrac{c_{p,2}}{r_2}\right) \le Q \\ \text{and } \dfrac{U_1r_1^2 r_2 + \left(r_2^2 U_2 + (-2U_1c_{p,1} - 2U_2c_{p,2})r_2 + U_2c_{p,2}^2\right)r_1 + U_1c_{p,1}^2 r_2}{4r_1 r_2} - K > \dfrac{U_1r_1^2 r_2 + (U_2r_2^2 + (-2U_1c_{m,1} - 2U_2c_{m,2})r_2 + U_2c_{m,2}^2)r_1 + U_1c_{m,1}^2 r_2}{4r_1 r_2} \\ \left(-\dfrac{U_1(c_{m,1} - r_1)}{2r_1}, -\dfrac{(c_{m,2} - r_2)U_2}{2r_2}, 0\right), & \text{otherwise} \end{cases}$$

In the first case of Example 5, 3D printing is more profitable; therefore, 3D printing technology is adopted (i.e., $v^* = 1$). However, the capacity of 3D printing is limited. Therefore, retailer's optimal order quantity is a function of 3D printing capacity (i.e., $q_i^* = \dfrac{U_i(2Qr_{i^-} - U_{i^-}c_{p,i} + U_{i^-}c_{p,i^-} + U_{i^-}r_i - U_{i^-}r_{i^-})}{2(U_i r_{i^-} + U_{i^-}r_i)}$). We will refer to this case as *"Case 1"* in the remainder of the paper. In the second case of Example 5, 3D printing is still more profitable, therefore 3D printing technology is still adopted (i.e., $v^* = 1$). Furthermore, the capacity of 3D printing is not limited. So, the retailer's optimal order quantity is very similar to that of the Newsvendor (i.e., $q_i^* = \dfrac{U_i}{2}\left(1 - \dfrac{c_{p,i}}{r_i}\right)$). We will refer to this case as *"Case 2"* in the remainder of the paper. In the third case of Example 5, 3D printing is not even profitable. Therefore, the problem simplifies to the Benchmark case. We will refer to this case as *"Case 3"* in the remainder of the paper. Finally, In the Numerical Study section, we examine the impact of various parameters on the equilibrium.

Figure 1 shows the equilibrium at $(U_1, U_2, C_{p,1}, C_{p,2}, Q, K, C_{m,1}, C_{m,2}) = (100, 150, 10, 20, 90, 400, 15, 30)$ with respect to $r_1$, $r_2$. One can observe that when $r_i$ are just a bit larger than $\max(C_{p,i}, C_{m,i})$, then traditional ways of manufacturing are opted for at the equilibrium. As either of $r_1$ or $r_2$ increases, then 3D printing will be adopted. When $r_i$ is not very large, the retailers optimal order quantity (i.e., $q_i$) is not very large because $q_i^* = \dfrac{U_i}{2}\left(1 - \dfrac{c_{p,i}}{r_i}\right)$. Therefore, 3D printing can be utilized with *no* capacity constraint. However, when $r_i$ is a lot larger than $\max(C_{p,i}, C_{m,i})$, then the retailers optimal order quantity (i.e., $q_i$) is going to be large at the equilibrium because $q_i^* = \dfrac{U_i}{2}\left(1 - \dfrac{c_{p,i}}{r_i}\right)$. As the result, 3D printing will be utilized at its capacity. In the numerical study section, we explore this topic in greater depth.

**Figure 1**- a) Manufacturer's Profit VS $r_1$ and $r_2$. b) Retailers' Profit VS $r_1$ and $r_2$

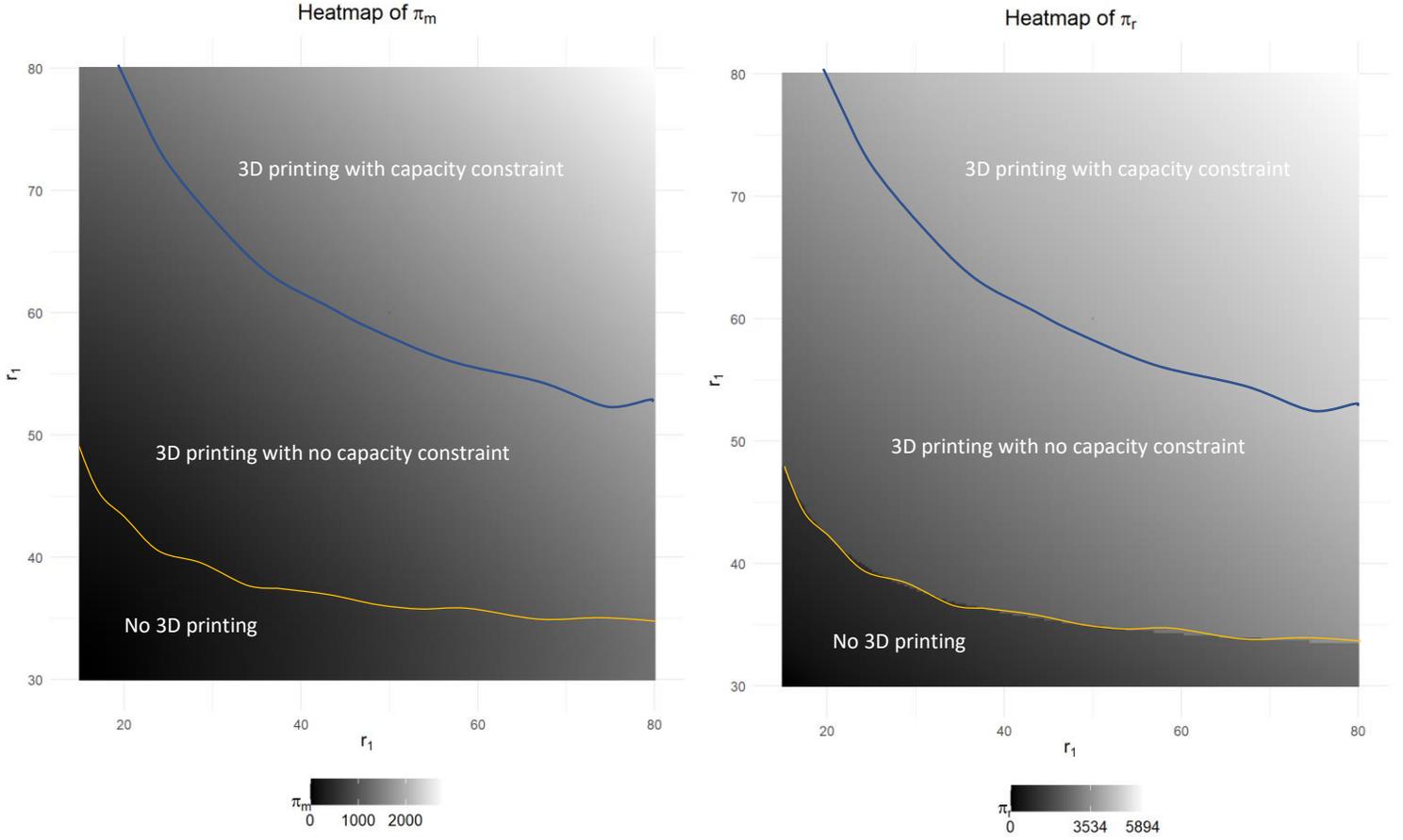

### 4.4. Investing in 3D Printing

In this section, we study under what economic conditions investing in 3D printing is beneficial for the supply chains, and our focus in this section is on the results for the Uniform distribution. In the numerical study section, we further investigate this question for other demand distributions.

Recall that 1) if 3D printing is adopted at capacity, $\pi_M^{1Q} = \frac{\left((c_{p,1} - c_{p,2} - r_1 + r_2)^2 U_2 - 4Qr_2(c_{p,1} - r_1)\right)U_1 - 4Qr_1\left((c_{p,2} - r_2)U_2 + Qr_2\right)}{4U_1 r_2 + 4U_2 r_1} - K$ (i.e., Case 1), 2) if 3D printing is adopted with no capacity constraint, $\pi_M^1 = \frac{U_1 r_1^2 r_2 + (r_2^2 U_2 + (-2U_1 c_{p,1} - 2U_2 c_{p,2})r_2 + U_2 c_{p,2}^2)r_1 + U_1 c_{p,1}^2 r_2}{4r_1 r_2} - K$ (i.e., Case 2), and 3) if 3D printing technology is not adopted the manufacturer's profit is $\pi_M^0 = \frac{U_1 r_1^2 r_2 + (U_2 r_2^2 + (-2U_1 c_{m,1} - 2U_2 c_{m,2})r_2 + U_2 c_{m,2}^2)r_1 + U_1 c_{m,1}^2 r_2}{4r_1 r_2}$ (i.e., Case 3). The following lemma, given the adoption of 3D printing technology, compares the incapacitated and capacitated problem.

**Lemma 3**. *For a uniform distribution, the unconstrained problem is always more profitable:*

$$\pi_M^1 - \pi_M^{1Q} =$$

$$\frac{U_1 r_1^2 r_2 + \left(r_2^2 U_2 + (-2U_1 c_{p,1} - 2U_2 c_{p,2})r_2 + U_2 c_{p,2}^2\right)r_1 + U_1 c_{p,1}^2 r_2}{4 r_1 r_2} - \frac{\left((c_{p,1} - c_{p,2} - r_1 + r_2)^2 U_2 - 4Q r_2 (c_{p,1} - r_1)\right) U_1 - 4Q r_1 \left((c_{p,2} - r_2)U_2 + Q r_2\right)}{4 U_1 r_2 + 4 U_2 r_1}$$

$$= \frac{\left(\left(\left(Q - \frac{U_1}{2} - \frac{U_2}{2}\right)r_2 + \frac{U_2 c_{p,2}}{2}\right)r_1 + \frac{U_1 c_{p,1} r_2}{2}\right)^2}{(U_1 r_2 + U_2 r_1) r_1 r_2} > 0$$

This lemma implies that, if 3D printing technology is adopted, a 3D printer whose capacity is more than the sum of retailer's Newsvendor solutions (i.e., $Q > q_{nv}(w_1) + q_{nv}(w_2)$) results in higher profits. So, given the opportunity, investing in 3D printing technology with higher capacity could be more profitable.

Next, we compare two cases where 3D printing is adopted with no capacity constraint and the case where 3D printing is not adopted.

**Lemma 4**. *For a uniform distribution, the condition under which adopting 3D printing technology is more profitable than using traditional ways of manufacturing is as follows.*

$$2K \leq U_1(c_{m,1} - c_{p,1})\left(1 - \frac{c_{m,1} + c_{p,1}}{2 r_1}\right) + U_2(c_{m,2} - c_{p,2})\left(1 - \frac{c_{m,2} + c_{p,1}}{2 r_2}\right)$$

This condition implies if the cost of adopting 3D printing technology is below a certain threshold, then it is beneficial to adopt 3D prating; otherwise, it is not.

Figure 2 summarizes Lemma 3 and Lemma 4 for the case where $(U_1, U_2, c_{m,1}, c_{m,2}, c_{p,1}, r_1, r_2, K, Q) = (10,15,5,10,1,10,20,10,8)$. It shows the profits in the above 3 cases versus $c_{p,2}$. The horizontal line reparents the manufacturer's profit when 3D printing is not adopted (i.e., Case 3). The blue and red curved lines represent the manufacturer's profit when 3D printing is adopted without capacity constraint (i.e., Case 2) and with capacity constraint (i.e., Case 1), respectively. One can observe that when $c_{p,2}$ is large (e.g., $c_{p,2} = 20$), then 3D printing is expensive, and this technology is not opted for in the equilibrium (i.e., the horizontal line is above the two curved lines). If 3D printing variable cost is below a specific threshold (i.e., $c_{p,2} \leq 15$ in Figure 2), however, then 3D printing is adopted. One can observe that in this region, the curved lines are above the horizontal line. The vertical line at $c_{p,2} = 10.7$ represents the capacity constraint boundary. That is, to the right of this line (i.e., $c_{p,2} \geq 10.7$), 3D printing variable cost is only moderately cheap. As a result, while 3D printing technology can be utilized, its moderately low per-unit cost limits the manufacturer's ability to set a very low wholesale price. Consequently, the retailer is unlikely to order a large quantity of products, meaning 3D printing's full capacity will not be reached. Finally, if the variable cost of 3D printing is sufficiently low (i.e., $c_{p,2} \leq 10.7$), a substantial number of products should be produced using 3D printing, leading to a scenario where the technology reaches its production capacity.

**Figure 2**-Manufacturer's profit for different cases VS $c_{p,2}$

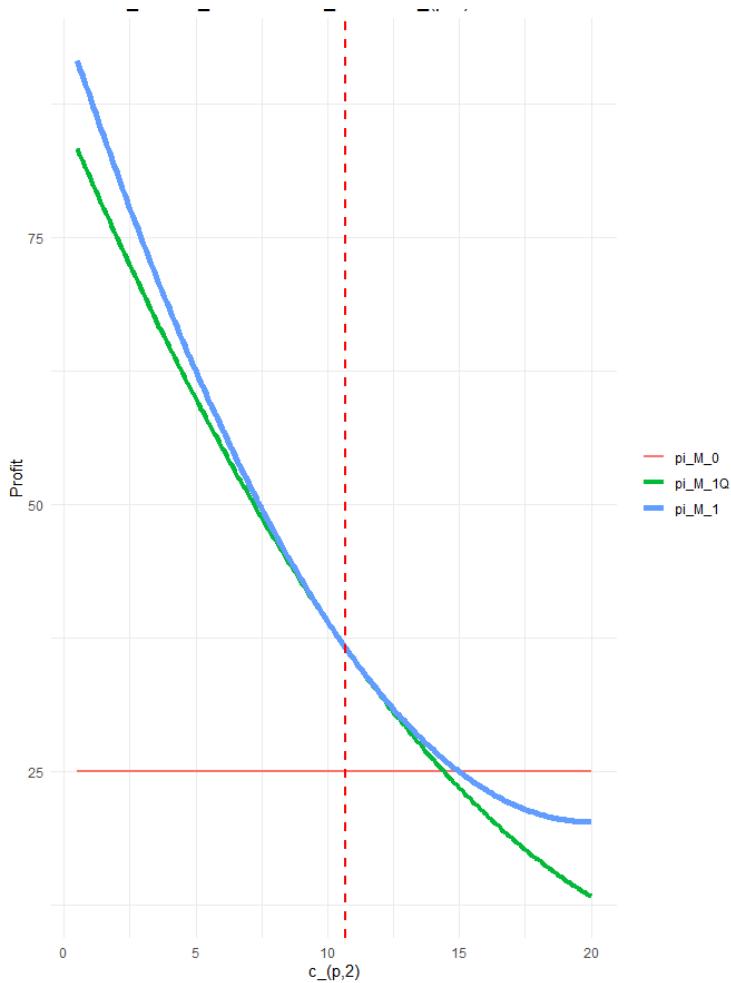

## 5. Manufacturer May Purchase a 3D Printer—Multiple Products

In this section, we investigate a scenario where the manufacturer uses traditional manufacturing methods to produce multiple distinct, $n$, products. The unit production cost for product $i$ is $c_{m,i}$. The manufacturer sells the products to the retailer at the wholesale price of $w_i; i \in \{1,2,\dots,n\}$. The retailer, in turn, faces stochastic demand for product $i$ with density functions of $f_i(.)$ and cumulative distribution functions of $F_i(.)$. The retailer sells the products to the final customer at the price of $r_i$. In this case, following our benchmark case, the retailer's problem is:

$$\pi_R^b(w) = \max_{q_i \geq 0} \sum_{i=1}^n E[r_i \min\{q_i, D_i\}] - w_i q_i$$

and the optimal solution to the retailer's problem is $q_{nv}(w_i) = F_i^{-1}\left(1 - \frac{w_i}{r_i}\right)$.

Next, similar to the previous section, in presence of 3D printing, the manufacturer has the option of purchasing a 3D printer with a fixed cost $K > 0$ and unit printing costs $c_{p,i}; i \in \{1,2,\dots,n\}$, to replace her traditional manufacturing. Here, again, we focus on cases where 3D printing replaces tradition manufacturing techniques.

Next, we let $v \in \{0,1\}$ denote the manufacturer's binary decision whether to adopt 3D printing technology. If 3D printing is adopted ($v = 1$), there will be no traditional manufacturing. Similarly, $Q$ is the capacity of 3D printing. Alternatively, if 3D printing is not adopted ($v = 0$), then only traditional manufacturing would be possible. Therefore, in presence of 3D printing, the manufacturer's profit-maximizing problem is

$$\pi_M^{3D} = \max_{w_i \geq 0, v \in \{0,1\}} \sum_{i=1}^n \left((w_i - v c_{p,i} - (1-v) c_{m,i}) q_{nv}(w_i)\right) - Kv \tag{10}$$

$$s.t. \quad \sum_{i=1}^n q_{nv}(w_i) \leq (1-v)M + vQ,$$

### 5.1. 3D Printing Is Not Adopted

In this section, we study the scenario where 3D printing is not adopted (i.e., $v = 0$), in which case, the problem in (10) simplifies to

$$\pi_M^{3D} = \max_{w_i \geq 0} \sum_{i=1}^n \left((w_i - c_{m,i}) q_{nv}(w_i)\right) \tag{11}$$

$$s.t. \quad \sum_{i=1}^n q_{nv}(w_i) \leq M,$$

Since $M$ is a large number, then technically, there is no constraint in (11). Therefore, the optimal solution can be achieved using the following proposition.

**Proposition 5**. *The unique optimal solution to the manufacturer's Problem (11) is*
$q_i^0 = q_b(c_{m,i}); i \in \{1,2,\dots,n\}$

$$w_i^0 = r_i\left(1 - F_i(q_i^0)\right); i \in \{1,2,\ldots,n\}$$

$$\pi_M^0 = \sum_{i=1}^{n}\left((w_i^0 - c_{m,i})q_b(c_{m,i})\right)$$

As one may observe, because in this case 3D printing is not adopted, the optimal order quantities and the optimal wholesale prices are similar to that of the benchmark (i.e., Proposition 1).

### 5.2. 3D Printing Is Adopted

In this section, we study the scenario where 3D printing is adopted (i.e., $v = 1$). In this case, the problem in (10) simplifies to

$$\pi_M^{3D} = \max_{w_i \geq 0} \sum_{i=1}^{n}\left((w_i - c_{p,i})q_{nv}(w_i)\right) - K \qquad (12)$$

$$s.t. \quad \sum_{i=1}^{n} q_{nv}(w_i) \leq Q,$$

Given the constraint in (12), we divide the solution space into two regions.

**Region 1.** $\sum_{i=1}^{n} q_{nv}(w_i) \leq Q$. In this case, the constraint in (12) is redundant, and there is an interior solution to the problem in (12). This case represents a scenario where 3D printing is *not* utilized at its maximum capacity. The following lemma maximizes manufacture's profit in this region.

**Lemma 5**. *The unique optimal solution to the manufacturer's Problem (12) in Region 1 is*
$$q_i^1 \triangleq q_b(c_{p,1}) = F_i^{-1}\left(1 - \frac{c_{p,i}}{r_i}\right)$$

$$w_i^1 \triangleq r_i\left(1 - F_i(q_i^1)\right)$$

$$\pi_M^{3D} \triangleq \sum_{i=1}^{n}\left((w_i^1 - c_{p,i})q_i^1\right) - K$$

**Region 2.** $\sum_{i=1}^{n} q_{nv}(w_i) \geq Q.$ In this case, the constraint in (12) is binding, and there is a corner solution to the problem in (12). This case represents a scenario where 3D printing *is* utilized at its maximum capacity. Therefore, $\sum_{i=1}^{n} q_{nv}(w_i) = Q$. The LaGrangian function of the problem in (12) is

$$L = \sum_{i=1}^{n}\left((w_i - c_{p,i})q_{nv}(w_i)\right) - K - \lambda\left(\sum_{i=1}^{n} q_{nv}(w_i) - Q\right)$$

We know, $w_i = r_i(1 - F_i(q_i))$. This simplifies the LaGrangian function to

$$L = \sum_{i=1}^{n}\left((r_i(1-F_i(q_i))-c_{p,i})q_i\right) - K - \lambda\left(\sum_{i=1}^{n} q_i - Q\right) \tag{13}$$

The optimal solution to Problem (13) is derived in the following Lemma.

**Lemma 6**. *The unique optimal solution to the manufacturer's Problem (13) in Region 2 is the unique solution to the following set of $(n+1)$ equations.*

$$\begin{cases} r_i(1-F_i(q_i)) - r_i q_i f_i(q_i) = c_{p,i} + \lambda; i \in \{1,2,\dots,n\} \\ \sum_{i=1}^{n} q_i = Q \end{cases}$$

Finally, the equilibrium for the Problem in (10) is derived in the following proposition.

**Proposition 6.** the equilibrium for the Problem in (10) is

$$(\vec{q_i^*}, v^*) = \begin{cases} (q_1^{1Q}, q_2^{1Q}, \dots, q_n^{1Q}, 1), & \sum_{i}^{n} q_i^1 > Q \text{ and } \pi_M^{1Q} > \pi_M^0 \\ (q_1^1, q_2^1, \dots, q_n^1, 1), & \sum_{i}^{n} q_i^1 \leq Q \text{ and } \pi_M^1 > \pi_M^0 \\ (q_1^0, q_2^0, \dots, q_n^0, 0), & otherwise \end{cases}$$

$$w_i^* = r_i(1 - F_i(q_i^*))$$

The first part of the equilibrium in Proposition 6 represents a scenario where 3D printing operates at full capacity (i.e., Case 1). The second part corresponds to a scenario where 3D printing is adopted without being limited by capacity constraints (i.e., Case 2). The third part reflects a scenario in which 3D printing is not adopted (i.e., Case 3).

## 6. Numerical Study

In this section, we conduct a numerical analysis of the equilibrium to provide deeper insights into optimal strategies for supply chain owners in presence of 3D printing technology. The parameters used for this study are listed in Table 1 and Table 2. Here, we focus on novel findings and avoid discussing results that are straightforward, such as the effect of $c_{m,i}$ on the equilibrium. We borrow the model in this section from Golchin and Rekabdar (2024).

### 6.1. Manufacturer May Purchase a 3D Printer—Two Products

In this section, we study 1) the effect of 3D printing variable costs (i.e., $c_{p,i}$), and 2) the effect of average demand (i.e., $U_i$) on the equilibrium. Specifically, we use the following problem parameters.

**Table 1**- Problem parameters for the numerical results for the two products scenarios.

| U1 | U2 | r1 | r2 | $c_{m,1}$ | $c_{m,2}$ | $c_{p,1}$ | $c_{p,2}$ | K | Q |
|---|---|---|---|---|---|---|---|---|---|
| 100 | $\{\frac{1}{2}U1, U1, 2U1\}$ = {50,100,200} | 50 | $\{\frac{1}{2}r1, r1, 2r1\}$ = {25,50,100} | $\{1, \ldots, r_1\}$ | $\{1, \ldots, r_2\}$ | $\{1, \ldots, r_1\}$ | $\{1, \ldots, r_2\}$ | {0,100,…,600} | {0,10,…,200} |

As there is more and more research on 3D printing, one can expect 3D printing variable cost to decrease. That is why in Observation 1, we specifically study how changes in $c_{p,i}$ affects the equilibrium.

**Observation 1.** As 3D printing variable cost (i.e., $c_{p,i}$) decreases, supply chain owners should opt for utilizing 3D printing technology. 3D printing can be adopted even if 3D printing is more expensive than traditional ways of manufacturing (i.e., $c_{p,i} > c_{m,i}$).

Figure 3 illustrates Observation 1, where three cases are presented: Case 1 (3D printing with capacity constraints), Case 2 (3D printing without capacity constraints), and Case 3 (no 3D printing). Figure 3.b overlays the condition $c_{p,i} = c_{m,i}$. Notably, when $c_{p,1}$ and $c_{p,2}$ significantly exceed $c_{m,1}$ and $c_{m,2}$ and approach $r_1$ and $r_2$, respectively, adopting 3D printing technology is suboptimal. Conversely, as $c_{p,i}$ decreases and diverges from $r_i$, 3D printing becomes viable. However, due to moderate $c_{p,i}$ values, retailers refrain from ordering large amounts, avoiding capacity saturation. Further decreases in $c_{p,i}$ lead to sustained 3D printing adoption, but reduced costs prompt retailers to increase orders, rendering 3D printer capacity a binding constraint.

Additionally, Figure 3.a displays the manufacturer's profit. A key observation is that, in the absence of 3D printing (Case 3), manufacturer's profit remains invariant with respect to $c_{p,i}$. Conversely, when 3D printing is employed, profit increases as $c_{p,i}$ decreases.

Observation 1 challenges existing literature (e.g., Arbabian, 2022), which posits that 3D printing adoption requires costs to be lower than traditional manufacturing methods (i.e., $c_{p,i} < c_{m,i}$). Our analysis reveals that this condition is not a necessary prerequisite for 3D printing adoption. Specifically, even when $c_{p,2} > c_{m,2} = 30$, 3D printing can still be viable if $c_{p,1}$ falls within certain ranges; specifically, if $c_{p,1} < 2$, 3D printing is adopted at full capacity, while $2 \leq c_{p,1} < 15$ yields 3D printing adoption without capacity constraints.

This phenomenon arises because 3D printing's benefits for one product can offset losses incurred by another. Figure 3.b illustrates this point. Contrary to Arbabian's (2022) assertion that 3D printing should not be adopted when $c_{p,1} = 20 > c_{m,1} = 15$, our analysis shows that if $c_{p,2} < 27$, 3D printing benefits product 2, outweighing the disadvantages for product 1. Consequently, 3D printing is adopted at equilibrium.

**Figure 3** – a) Manufacturer's profit. b) Equilibrium cases. $(U_1, U_2, c_{m,1}, c_{m,2}, r_1, r_2, K, Q) = (100, 150, 15, 30, 50, 100, 0, 100)$

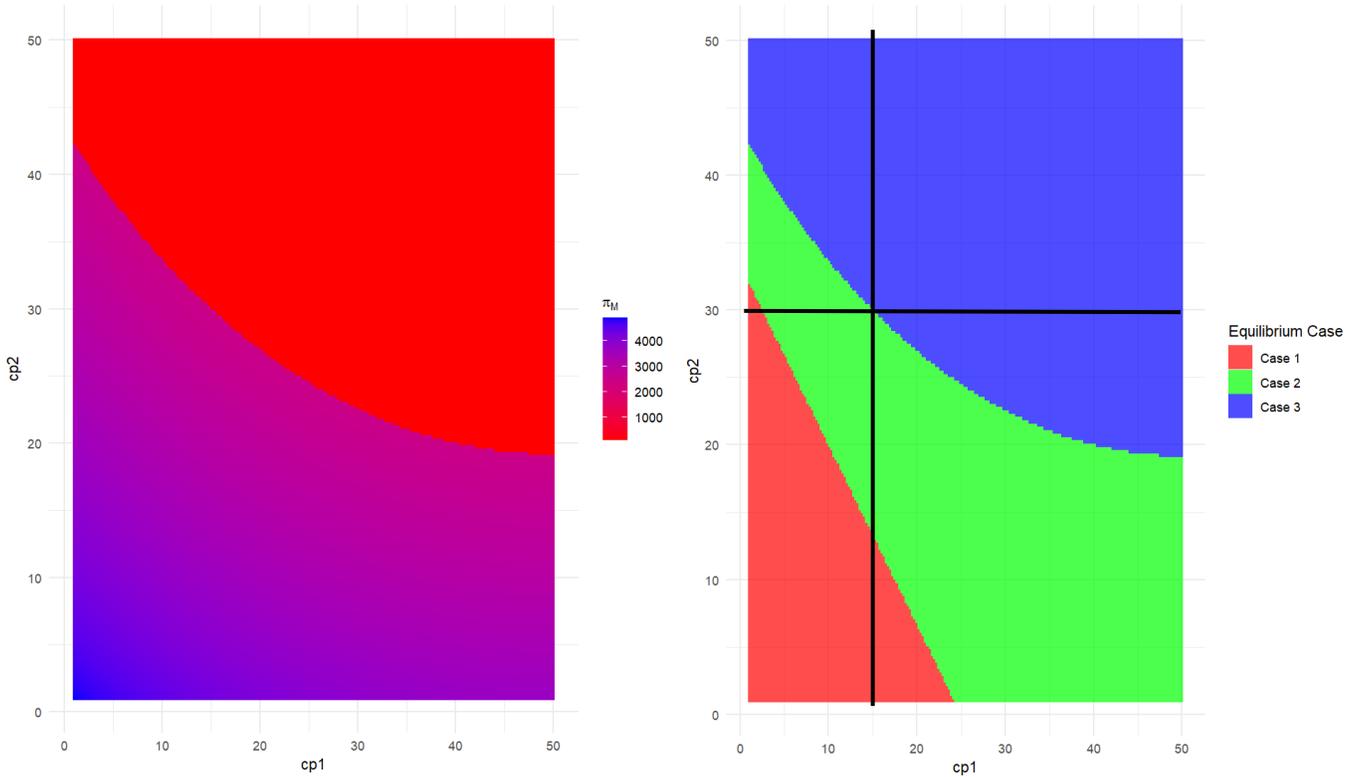

Note that on Figure 3, $K = 0$. We tried different $K$ values. When $K$ increases, a comparative analysis reveals that the region advocating 3D printing adoption contracts significantly. This shrinkage underscores an essential insight: as research and development in 3D printing advance, the associated adoption costs decrease. Consequently, this downward cost trajectory should incentivize an increasing number of companies to integrate 3D printing technology into their operations.

**Observation 2.** As average demand increases, supply chain owners should consider shifting from traditional manufacturing to 3D printing. However, when demand surpasses a certain threshold, traditional manufacturing becomes more advantageous.

Observation 2, illustrated in Figure 4, posits that when the average demand for products 1 and 2 is relatively low, supply chain managers should continue utilizing traditional manufacturing methods, as the adoption of 3D printing technology incurs significant costs (i.e., $K > 0$). This is evident in Figure 4.b, where for small values of $U_1$ and $U_2$, Case 3 (i.e., no 3D printing) emerges as the optimal equilibrium. As $U_i$ increases, the equilibrium shifts from traditional manufacturing to 3D printing, driven by rising average sales, which in turn elevate the supply chain's expected profit and justify the cost of adopting 3D printing technology. In this scenario, because the average demand remains moderate, 3D printing operates without capacity constraints. Case 2 of Figure 4 illustrates this case. As $U_i$ continues to grow, 3D printing is still preferred, but the system begins to experience capacity limitations, with the technology's bottleneck becoming apparent (Case 3 of Figure 4 illustrates this case). The counterintuitive finding arises when $U_i$ surpasses a critical threshold (the top right corner of Figure 4b, where case 3 takes over Case 1): beyond this point, 3D printing should not be utilized. In this high-demand scenario, the limited capacity of 3D printers restricts potential sales, leading to revenue loss if 3D printing is adopted. Consequently, the equilibrium favors traditional manufacturing methods, which are better suited to meet the increased demand efficiently.

**Figure 4** – a) Manufacturer's profit VS $U_1$ and $U_2$. b) Equilibrium cases. $(c_{p,1}, c_{p,2}, c_{m,1}, c_{m,2}, r_1, r_2, K, Q) = (10, 20, 15, 30, 50, 100, 400, 100)$

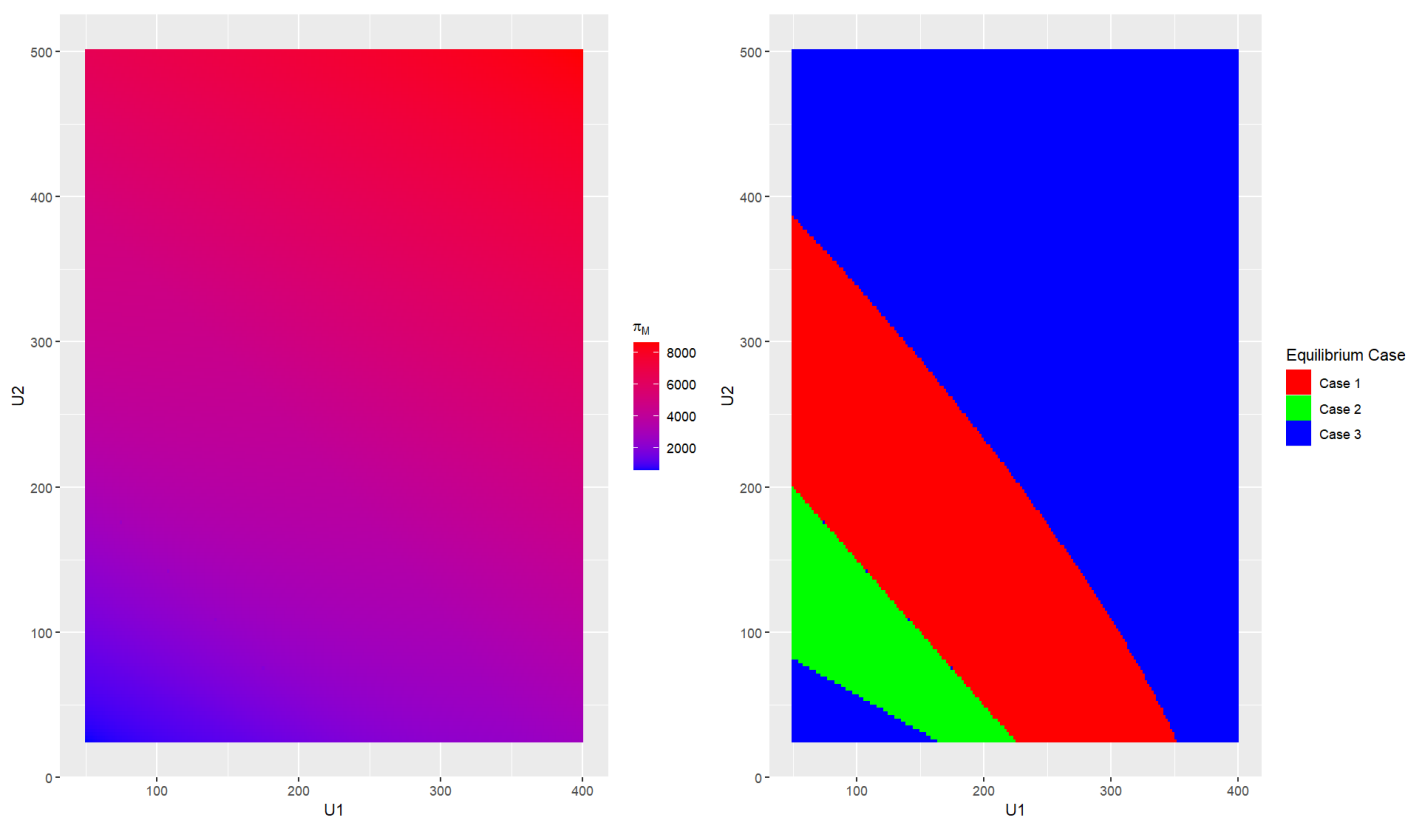

**6.2. Manufacturer May Purchase a 3D Printer—Multiple Products**

In this section, we focus on 1) uniform demand distributions and 2) a supply chain that produces and sells 3 distinctive products. By comparing the results in this section to the results in Section 6.1., we study the effect of number of products on the equilibrium. Furthermore, on top of the problem parameters used in the previous section, we use the following problem parameters.

**Table 2**- Problem parameters for the numerical results for the three-product scenario

| $U_3$ | $r_3$ | $c_{m,3}$ | $c_{p,3}$ |
|---|---|---|---|
| {50,100,200,400} | {25,100,150, 200,300} | $\{1, \dots, r_3\}$ | $\{1, \dots, r_3\}$ |

Following the results in Proposition 6, one can derive the following.

$$q_1^{1Q} = \frac{U_1(2Qr_2r_3 + U_2c_{p,2}r_3 - U_2c_{p,1}r_3 - U_2r_2r_3 + U_2r_1r_3 - U_3c_{p,1}r_2 + U_3c_{p,3}r_2 + U_3r_1r_2 - U_3r_2r_3)}{2(U_1r_2r_3 + U_2r_1r_3 + U_3r_1r_2)}$$

$$q_2^{1Q} = \frac{U_2(2Qr_1r_3 + U_1c_{p,1}r_3 - U_1c_{p,2}r_3 - U_1r_1r_3 + U_1r_2r_3 - U_3c_{p,2}r_1 + U_3c_{p,3}r_1 + U_3r_1r_2 - U_3r_1r_3)}{2(U_1r_2r_3 + U_2r_1r_3 + U_3r_1r_2)}$$

$$q_3^{1Q} = \frac{U_3(2Qr_1r_2 + U_1c_{p,1}r_2 - U_1c_{p,3}r_2 - U_1r_1r_2 + U_1r_2r_3 - U_2c_{p,3}r_1 + U_2c_{p,2}r_1 + U_2r_1r_3 - U_2r_1r_2)}{2(U_1r_2r_3 + U_2r_1r_3 + U_3r_1r_2)},$$

$$q_i^1 = U_i\left(1 - \frac{c_{p,i}}{r_i}\right); i \in \{1,2,3\},$$

$$q_i^0 = U_i\left(1 - \frac{c_{m,i}}{r_i}\right); i \in \{1,2,3\},$$

$$w_i^* = r_i\left(1 - \frac{q_i^*}{U_i}\right).$$

**Observation 3.** As the number of products in a supply chain grows, the operational regions where 3D printing is employed expand.

Figure 5 illustrates Observation 3, with Figure 5.b depicting the equilibrium behavior relative to $c_{p,3}$. In this context, Case 1 represents the scenario where 3D printing is used at full capacity, Case 2 represents the scenario where 3D printing is used without capacity constraints, and Case 3 refers to the scenario where 3D printing is not utilized at all. The parameters used are intentionally similar to those in Figure 3, where $c_{m,1} = 15$ and $c_{m,2} = 30$. Additionally, we set $c_{p,1} = 16 > c_{m,1}$ and $c_{p,2} = 31 > c_{m,2}$, which corresponds to the region in Figure 3 where 3D printing is *NOT* adopted (i.e., Case 3). However, in Figure 5, 3D printing is *utilized* when $c_{p,3} \leq 45$. This outcome arises because, although 3D printing may not be advantageous for products 1 and 2, its benefits for product 3 are substantial enough to compensate for the potential losses incurred for products 1 and 2.

Figure 5.b illustrates the manufacturer's profits as a function of $c_{p,3}$. As expected, the profit decreases as $c_{p,3}$ increases. On this graph, the lines indicating where the equilibrium transitions between Case 1, Case 2, and Case 3 have been superimposed. Specifically, when $c_{p,3} < 24$, 3D printing is used with a capacity constraint. For $24 \leq c_{p,3} < 44$, 3D printing is utilized without a capacity constraint. Finally, when $c_{p,3} \geq 44$, 3D printing is no longer employed.

**Figure 5** – a) Manufacturer's profit VS $c_{p,3}$. b) Equilibrium cases. $(c_{m,3}, r_1, r_2, r_3, K, Q) = (45, 50, 100, 150, 0, 170)$

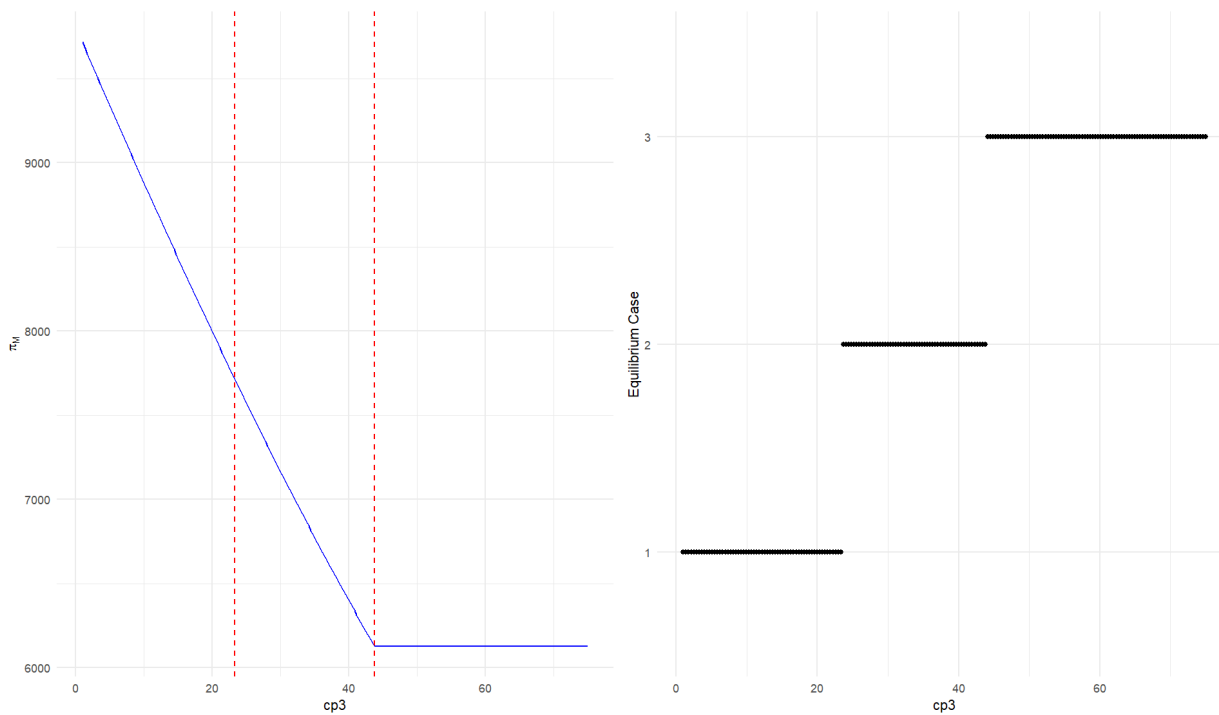

Observation 3 underscores the critical role of integrating 3D printing within supply chains characterized by product customization. Specifically, as indicated in this observation, there is a positive correlation between the number of distinct products produced by a supply chain and the extent of 3D printing adoption. In scenarios where supply chains offer a high variety of products (tailored to individual customer preferences), the expectation is that the implementation of 3D printing technologies will occur more rapidly compared to contexts where customization is minimal or absent. This alignment with contemporary literature highlights the potential of additive manufacturing to enhance responsiveness and flexibility in meeting diverse consumer demands (García-Alonso et al., 2022; Top et al., 2023).

## 7. Conclusion

This study analyzes the integration of 3D printing into decentralized supply chains, focusing on the viability of adoption despite the fixed cost of investing in 3D printing technology. Utilizing a Stackelberg game framework, we examine two primary scenarios: (1) a two-product production model and (2) a multi-product model, each constrained by capacity. In the first scenario, we derive a sufficient condition to identify equilibrium for a generic demand distribution, and for a uniform distribution, we fully characterize the equilibrium. Additionally, we determine a threshold for the fixed costs of adopting 3D printing, below which 3D printing becomes the more profitable choice. For the multi-product case, we also establish a sufficient condition to derive the equilibrium.

Our numerical results reveal three key observations regarding the adoption of 3D printing in supply chains. First, as the number of products in a supply chain increases, the operational regions where 3D printing is viable expand, even if 3D printing is more expensive than traditional manufacturing methods. This is because the benefits of 3D printing for certain products can outweigh the costs associated with others. Second, as the variable costs of 3D printing decrease, supply chains are more likely to adopt the technology, though capacity constraints become a critical factor. When the cost of 3D printing decreases significantly, retailers increase orders, pushing 3D printing capacity to its limits. Lastly, 3D printing becomes more favorable as demand increases, but only up to a certain threshold. Beyond that, traditional manufacturing methods are more efficient, as limited 3D printing capacity restricts the ability to meet high demand.

Future research could build on these findings by investigating how technological advancements in 3D printing—such as reduced setup times and cost efficiencies—might shift adoption thresholds even in high-demand scenarios. Additionally, examining the impact of 3D printing on environmental sustainability and localized production strategies could provide a broader perspective on its implications for supply chain resilience. By understanding these evolving dynamics, supply chain managers can better assess when and how to deploy 3D printing to achieve strategic goals within complex market environments.

**Data Availability Statement**

The authors confirm that the data supporting the findings of this study are available within the paper.

**References**

1. Ahmed, R., Heese, H. S., and Kay, M. 2023. "Designing a Manufacturing Network with Additive Manufacturing Using Stochastic Optimization." *International Journal of Production Research* 61 (7): 2267–2287.
   https://doi.org/10.1080/00207543.2022.2056723
2. Arbabian, M., and Wagner, M. 2020. "The Impact of 3D Printing on Manufacturer–Retailer Supply Chains." *European Journal of Operational Research* 285 (2): 538–552.
   https://doi.org/10.1016/j.ejor.2020.01.063
3. Arbabian, M. 2022. "Supply Chain Coordination via Additive Manufacturing." *International Journal of Production Economics* 243: 108318.
   https://doi.org/10.1016/j.ijpe.2021.108318
4. Banciu, M, and Mirchandani, P. 2013. "Technical Note: New Results Concerning Probability Distributions with Increasing Generalized Failure Rates." *Operations Research* 61 (4): 925–931.
   https://doi.org/10.1287/opre.2013.1198
5. Bernstein, F. and Kök, A. G. 2009. "Dynamic Cost Reduction through Process Improvement in Assembly Networks." *Management Science* 55 (4): 552–567.
   https://doi.org/10.1287/mnsc.1080.0961
6. Cantini, A., Peron, M., De Carlo, F., and Sgarbossa, F. 2022. "A Decision Support System for Configuring Spare Parts Supply Chains Considering Different Manufacturing Technologies." *International Journal of Production Research 62(8)*: 3023–3043.
   https://doi.org/10.1080/00207543.2022.2041757
7. Cachon, G. 2003. "Supply Chain Coordination with Contracts." In *Handbooks in Operations Research and Management Science* 11: 227–339.
   https://doi.org/10.1016/S0927-0507(03)11006-7
8. Cachon, G., and Lariviere, M. 2005. "Supply Chain Coordination with Revenue-Sharing Contracts: Strengths and Limitations." *Management Science* 51 (1): 30–44.
   https://doi.org/10.1287/mnsc.1040.0215
9. Chen, L., Cui, Y., and Lee, H. 2021. "Retailing with 3D printing". *Productions and Operations Management*. 30 (7): 1986-2007
   https://doi.org/10.1111/poms.13367
10. Cohen, M., Lobel, R., and Perakis, G. 2016. "The Impact of Demand Uncertainty on Consumer Subsidies for Green Technology Adoption." *Management Science* 62 (5): 1235–1258.
    https://doi.org/10.1287/mnsc.2015.2173
11. Cui, W., Yang, Y., and Di, L. 2023. "Modeling and Optimization for Static-Dynamic Routing of a Vehicle with Additive Manufacturing Equipment." *International Journal of Production Economics* 257: 108756.
    https://doi.org/10.1016/j.ijpe.2022.108756
12. Demiralay, E., Razavi, M., Kucukkoc, I., and Peron, M. 2023. "An Environmental Decision Support System for Determining On-Site or Off-Site Additive Manufacturing of Spare Parts." *IFIP International Conference on Advances in Production Management Systems* 563–574.
13. Dong, L., Shi, D., and Zhang, F. 2017. "3D Printing vs. Traditional Flexible Technology: Implications for Operations Strategy." *Manufacturing & Service Operations Management* 19 (2): 281–294.
    https://doi.org/10.1287/msom.2017.0623

14. Emelogu, A., Marufuzzaman, M., Thompson, S. M., Shamsaei, N., and Bian, L. 2016. "Additive Manufacturing of Biomedical Implants: A Feasibility Assessment via Supply-Chain Cost Analysis." *Additive Manufacturing* 11: 97–113.
    https://doi.org/10.1016/j.addma.2016.04.006
15. Fine, C., and Freund, R. 1990. "Optimal Investment in Product-Flexible Manufacturing Capacity." *Management Science* 36 (4): 449–466.
16. Golchin, B., and Rekabdar, B., Anomaly detection in time series data using reinforcement learning, variational autoencoder, and active learning, Proceedings of the 2024 Conference on AI, Science, Engineering, and Technology (AIxSET), pp. 1-8, Laguna Hills, CA, September 30 – October 2, 2024.
17. Goyal, M., and Netessine, S. 2007. "Strategic Technology Choice and Capacity Investment under Demand Uncertainty." *Management Science* 53 (2): 192–207.
    https://doi.org/10.1287/mnsc.1060.0611
18. Gupta, S. 2008. "Research Note - Channel Structure with Knowledge Spillovers." *Marketing Science* 27 (2): 247–261.
    https://doi.org/10.1287/mksc.1070.0285
19. Ishii, A. 2004. "Cooperative Research and Development Between Vertically Related Firms with Spillovers." *International Journal of Industrial Organization* 22 (8): 1213–1235.
    https://doi.org/10.1016/j.ijindorg.2004.05.003
20. Ivanov, D. 2021. "Supply Chain Viability and the COVID-19 Pandemic: A Conceptual and Formal Generalisation of Four Major Adaptation Strategies." *International Journal of Production Research* 59 (12): 3535–3552.
    https://doi.org/10.1080/00207543.2021.1890852
21. Ivanov, D., Dolgui, A., and Sokolov, B. 2019. "The Impact of Digital Technology and Industry 4.0 on the Ripple Effect and Supply Chain Risk Analytics." *International Journal of Production Research* 57 (3): 829–846.
    https://doi.org/10.1080/00207543.2018.1488086
22. Ivanov, D., Dolgui, A., Blackhurst, J., and Choi, T. 2023. "Toward Supply Chain Viability Theory: From Lessons Learned Through COVID-19 Pandemic to Viable Ecosystems." International Journal of Production Research 61 (8): 2402–2415.
    https://doi.org/10.1080/00207543.2023.2177049
23. Jarrar, Q., Belkadi, F., Blanc, R., Kestaneci, K., and Bernard, A. 2023. "Knowledge Reuse for Decision Aid in Additive Manufacturing: Application on Cost Quotation Support." *International Journal of Production Research* 61 (20): 7027–7047.
    https://doi.org/10.1080/00207543.2022.2142861
24. *Kumar, M.,* Graham, G., Hennelly, P., and Srai, J. 2016. "How Will Smart City Production Systems Transform Supply Chain Design: A Product-Level Investigation." *International Journal of Production Research 54 (23): 7181–7192.*
    https://doi.org/10.1080/00207543.2016.1198057
25. *Kunovjanek, M.,* and Reiner, G. 2020. "How Will the Diffusion of Additive Manufacturing Impact the Raw Material Supply Chain Process?" *International Journal of Production Research* 58 (5): 1540–1554.
    https://doi.org/10.1080/00207543.2019.1661537
26. Kucukkoc, I. 2019. "MILP Models to Minimise Makespan in Additive Manufacturing Machine Scheduling Problems." *Computers & Operations Research* 105: 58–67.
    https://doi.org/10.1016/j.cor.2019.01.006

27. Lolli, F., Coruzzolo, A. M., Peron, M., and Sgarbossa, F. 2022. "Age-Based Preventive Maintenance with Multiple Printing Options." *International Journal of Production Economics* 243: 108339.
https://doi.org/10.1016/j.ijpe.2021.108339
28. Pasternack, B. 1985. "Optimal Pricing and Return Policies for Perishable Commodities." *Marketing Science* 31 (8): 166–176.
https://doi.org/10.1287/mksc.4.2.166
29. Peron, M., Sgarbossa, F., Ivanov, D., and Dolgui, A. 2022. "Impact of Additive Manufacturing on Supply Chain Resilience During COVID-19 Pandemic." *Supply Network Dynamics and Control* 20: 121–146.
https://doi.org/10.1007/978-3-031-09179-7_6
30. Rodríguez-Espíndola, O., Chowdhury, S., Beltagui, A., and Albores, P. 2020. "The Potential of Emergent Disruptive Technologies for Humanitarian Supply Chains: The Integration of Blockchain, Artificial Intelligence, and 3D Printing." *International Journal of Production Research* 58 (15): 4610–4630.
https://doi.org/10.1080/00207543.2020.1761565
31. Sgarbossa, F., Peron, M., Lolli, F., and Balugani, E. 2021. "Conventional or Additive Manufacturing for Spare Parts Management: An Extensive Comparison for Poisson Demand." International Journal of Production Economics 233: 107993.
https://doi.org/10.1016/j.ijpe.2020.107993
32. Strong, D., Kay, M., Conner, B., Wakefield, T., and Manogharan, G. 2019. "Hybrid Manufacturing—Locating AM Hubs Using a Two-Stage Facility Location Approach." *Additive Manufacturing* 25: 469–476.
https://doi.org/10.1016/j.addma.2018.11.027
33. Tareq, M. S., Rahman, T., Hossain, M., and Dorrington, P. 2021. "Additive Manufacturing and the COVID-19 Challenges: An In-Depth Study." *Journal of Manufacturing Systems* 60: 787–798.
https://doi.org/10.1016/j.jmsy.2020.12.021
34. Thomas-Seale, L. E. J., Kirkman-Brown, J. C., Attallah, M. M., Espino, D. M., and Shepherd, D. E. T. 2018. "The Barriers to the Progression of Additive Manufacture: Perspectives from UK Industry." *International Journal of Production Economics* 198: 104–118.
https://doi.org/10.1016/j.ijpe.2018.02.003
35. Top, N., Sahin, I., Mangla, S. K., Sezer, M. D., and Kazancoglu, Y. 2023. "Towards Sustainable Production for Transition to Additive Manufacturing: A Case Study in the Manufacturing Industry." International Journal of Production Research 61 (13): 4450–4471.
https://doi.org/10.1080/00207543.2022.2152895
36. Van Mieghem, J. 1998. "Investment Strategies for Flexible Resources." *Management Science* 44 (8): 1071–1078.
https://doi.org/10.1287/mnsc.44.8.1071

**Appendix**

*Proof of Lemma 2.*

$$\frac{\partial \pi_M^{3D}}{\partial q_1} = -r_1 f_1(q_1) q_1 + \left(r_1\left(1 - F_1(q_1)\right) - c_{p,1}\right) + r_2 f_2(Q - q_1)(Q - q_1) - \left(r_2\left(1 - F_2(Q - q_1)\right) - c_{p,2}\right)$$

$$\to \frac{\partial^2 \pi_M^{3D}}{\partial q_1^2} = -r_1 f_1(q_1) - r_2 f_2(Q - q_1) - r_2 f_2(Q - q_1) < 0$$

We define $q_1^{1Q}$ as the unique solution to the following equation.

$$r_1(1 - F_1(q_1))\left[1 - \frac{f_1(q_1) q_1}{1 - F_1(q_1)}\right] - r_2(1 - F_2(Q - q_1))\left[1 - \frac{f_2(Q - q_1)(Q - q_1)}{1 - F_2(Q - q_1)}\right] = c_{p,1} - c_{p,2}.$$

*Proof of Lemma 5.*

$$\frac{\partial L}{\partial q_i} = -r_i q_i f_i(q_i) + r_i(1 - F_i(q_i)) - c_{p,i} - \lambda$$

$$\to \frac{\partial L}{\partial \lambda} = Q - \sum_{i=1}^{n} q_i$$

$$\to \frac{\partial^2 L}{\partial q_i^2} = -2 r_i f_i(q_i)$$

$$\to \frac{\partial^2 L}{\partial q_i \partial \lambda} = -1$$

$$\to H = \begin{bmatrix} -2 r_1 f_1(q_1) & \cdots & -1 \\ \vdots & \ddots & \vdots \\ -1 & \cdots & -2 r_n f_n(q_n) \end{bmatrix}$$

A straightforward Analysis shows that all of H's principal minors of even order are positive and all of its principal minors of odd order are negative.

$$m_1 = -2 r_1 f_1(q_1) < 0$$

$$m_2 = \det \begin{bmatrix} -2 r_1 f_1(q_1) & 0 \\ 0 & -2 r_2 f_2(q_2) \end{bmatrix} > 0$$

$$m_3 = \det \begin{bmatrix} -2 r_1 f_1(q_1) & 0 & 0 \\ 0 & -2 r_2 f_2(q_2) & 0 \\ 0 & 0 & -2 r_3 f_3(q_3) \end{bmatrix} < 0$$

This trend continues because all the main diagonal elements of the above matrix are negative, and all other elements are zero. Now, the only thing that is left to check is the last principal minor because it includes the derivative of $\lambda$. Let us assume n=2. In this case,

$$m_3 = \det \begin{bmatrix} -2r_1f_1(q_1) & 0 & -1 \\ 0 & -2r_2f_2(q_2) & 0 \\ -1 & 0 & 0 \end{bmatrix} = -1 * \det \begin{bmatrix} -2r_1f_1(q_1) & 0 \\ 0 & -2r_2f_2(q_2) \end{bmatrix} < 0.$$

Or, if n=3, then

$$m_4 = \det \begin{bmatrix} -2r_1f_1(q_1) & 0 & 0 & -1 \\ 0 & -2r_2f_2(q_2) & 0 & 0 \\ 0 & 0 & -2r_3f_3(q_3) & 0 \\ -1 & 0 & 0 & 0 \end{bmatrix}$$

$$= -2r_1f_1(q_1) \det \begin{bmatrix} -2r_2f_2(q_2) & 0 & 0 \\ 0 & -2r_3f_3(q_3) & 0 \\ 0 & 0 & 0 \end{bmatrix} + \det \begin{bmatrix} -2r_1f_1(q_1) & 0 & 0 \\ 0 & -2r_2f_2(q_2) & 0 \\ 0 & 0 & -2r_3f_3(q_3) \end{bmatrix}$$

$$< 0$$

Therefore, the optimal solution to the problem in 2-3 is the unique solution to the following set of equations:

$$\begin{cases} r_i(1 - F_i(q_i)) - r_iq_if_i(q_i) = c_{p,i} + \lambda; i \in \{1,2,\dots,n\} \\ \sum_{i=1}^{n} q_i = Q \end{cases}$$